\documentclass{article}

\usepackage{amssymb, latexsym}

\usepackage{graphicx}

\usepackage{amsmath}
\usepackage{color}

\newtheorem{theorem}{Theorem}[section]

\newtheorem{corollary}[theorem]{Corollary}

\newtheorem{proposition}[theorem]{Proposition}

\newtheorem{remark}[theorem]{Remark}

\begin{document}
\title{\large\bf  $W$-entropy, super Perelman Ricci flows and $(K, m)$-Ricci solitons} 

\author{\ \ Songzi Li\thanks{Research partially supported by a Postdoctral Fellowship of the Beijing Normal University.} , Xiang-Dong Li
\thanks{Research supported by NSFC No. 11771430, and Key Laboratory RCSDS, CAS, No. 2008DP173182.} }

\maketitle

\begin{minipage}{120mm}
{\bf Abstract}.  In this paper, we prove the characterization  of  the 
 $(K, \infty)$-super Perelman Ricci flows by various functional inequalities and  gradient estimate for the heat semigroup generated by   the Witten Laplacian on manifolds equipped with time dependent metrics and potentials.  
 As a byproduct, we derive the Hamilton type dimension free 
   Harnack inequality  on manifolds with  $(K, \infty)$-super 
Perelman Ricci flows.  Based on a new second order  differential inequality  on  the Boltzmann-Shannon entropy for the heat equation of the Witten Laplacian, we introduce a new $W$-entropy quantity and prove its monotonicity 
for the heat equation of the  Witten Laplacian on complete Riemannian manifolds with the  $CD(K, \infty)$-condition and on 
compact manifolds  with  $(K, \infty)$-super Perelman  Ricci flows.  
Our results characterize  the $(K, \infty)$-Ricci solitons and  the $(K, \infty)$-Perelman Ricci flows. We also prove a second order differential entropy inequality on $(K, m)$-super Ricci flows, 
which can be used to characterize the $(K, m)$-Ricci solitons and the $(K, m)$-Ricci flows. Finally, we give a probabilistic interpretation of the $W$-entropy for the heat equation of the Witten Laplacian on 
manifolds with the $CD(K, m)$-condition. 

 \end{minipage}

\bigskip


\medskip

$(K, m)$-super Perelman Ricci flow, $(K, m)$-Perelman Ricci flows.


\section{Introduction}

In \cite{H1}, R. Hamilton introduced the Ricci flow to deform Riemannian metrics on a manifold $M^n$ by the evolution equation 
\begin{eqnarray}
\partial_t g=-2Ric_{g}.
\end{eqnarray}
The volume preserving 
 normalized Ricci flow equation on a closed manifold $M^n$ is given by 
 \begin{eqnarray}
\partial_t g=-2Ric_{g}-{2\over n}rg, 
\end{eqnarray}
where $r={1\over {\rm Vol}(M)} \int_M Rdv$ is the average of the scalar curvature $R$ on $(M, g)$. In the case of $3$-dimensional closed manifolds with a metric $g_0$ of positive Ricci curvature , Hamilton \cite{H1} proved that 
the unique solution of the normalized Ricci flow $g(t)$  with $g(0)=g_0$ exists on $[0, \infty)$ and  converge exponentially fast in every $C^k$-norm, $k\in \mathbb{N}$, to a metric of positive constant sectional curvature. 
As a consequence,  the Poincar\'e conjecture is proved on simply connected $3$-dimensional  closed Riemannian manifolds with positive Ricci curvature.

In \cite{P1},  Perelman gave a gradient flow  interpretation for the Ricci flow and proved two entropy monotonicity results along the Ricci flow. More precisely, let $M$ be a closed manifold, $n=dim M$, define 
\begin{eqnarray*}
\mathcal{F}(g, f)=\int_M (R+|\nabla f|^2)e^{-f}dv,
\end{eqnarray*}
where $g\in \mathcal{M}=\{{\rm Riemannian\ metric \ on}\ M\}$,
$f\in C^\infty(M)$.  Under the constraint condition
that the weighted volume measure 
\begin{eqnarray*}
d\mu=e^{-f}dv\end{eqnarray*}
is fixed, Perelman \cite{P1} proved that the gradient flow of $\mathcal{F}$ with respect to the standard $L^2$-metric on  $\mathcal{M}\times C^\infty(M)$ is given by the following modified Ricci flow for $g$ together with  the conjugate heat equation for $f$, i.e.,
\begin{eqnarray*}
\partial_t g&=&-2(Ric+\nabla^2 f),\\
\partial_t f&=&-\Delta f-R.
\end{eqnarray*}
Moreover,   Perelman introduced the
remarkable $W$-entropy for the Ricci flow as follows
\begin{eqnarray}
W(g, f, \tau)=\int_M \left[\tau(R+|\nabla
f|^2)+f-n\right]{e^{-f}\over
 (4\pi\tau)^{n/2}}dv,\label{entropy-1}
 \end{eqnarray}
 and proved 
 the following beautiful 
$W$-entropy formula 
\begin{eqnarray}
{d\over dt}W(g, f, \tau)=2 \int_M \tau\left|Ric+\nabla^2
f-{g\over 2\tau}\right|^2{e^{-f}\over (4\pi \tau)^{n/2}}dv\label{Entropy-P}
\end{eqnarray}
along the evolution
equation
\begin{eqnarray}
\partial_t g=-2Ric, \ \ \partial_t f=-\Delta f+|\nabla
f|^2-R+{n\over 2\tau}, \ \ \partial_t \tau=-1.\label{r-c}
\end{eqnarray}
In particular,  the $W$-entropy is monotonic  increasing in $t$
and the monotonicity is strict except that $(M, g(\tau), f(\tau))$
is a shrinking Ricci soliton, i.e.,
\begin{eqnarray}
Ric+\nabla^2 f={g\over 2\tau}.\label{SRS}
\end{eqnarray}
As an application,  Perelman \cite{P1} proved the no local
collapsing theorem, which ``removes the major stumbling block in
Hamilton's approach to geometrization'' and plays an important role in the final resolution of the
Poincar\'e conjecture and Thurston's geometrization conjecture.

To better describe  the motivation and our results,  we need to introduce some definitions and notations. Let $(M, g)$ be a complete Riemannian manifold, $\phi\in C^2(M)$, and $d\mu=e^{-\phi}dv$, where $dv$ is the Riemannian volume measure on $(M, g)$. The Witten Laplacian, called also the weighted Laplacian, and denoted by 
\begin{eqnarray*}
L =\Delta -\nabla \phi\cdot\nabla \label{WL}
\end{eqnarray*}
is a self-adjoint and non-positive operator on $L^2(M, \mu)$. For all $u, v\in C^\infty_0(M)$, the following integration by parts formula holds
\begin{eqnarray*}
\int_M \langle \nabla u, \nabla v\rangle d\mu=-\int_M Lu vd\mu=-\int_M uLvd\mu.
\end{eqnarray*}
%
In \cite{BE}, Bakry and Emery proved that for all $u\in C_0^\infty(M)$,
\begin{eqnarray}
L|\nabla u|^2-2\langle \nabla u, \nabla L u\rangle=2|\nabla^2
u|^2+2Ric(L)(\nabla u, \nabla u), \label{BWF}
\end{eqnarray}
where
$$Ric(L)=Ric+\nabla^2\phi,$$
which appeared in the Ricci soliton equation $(\ref{SRS})$ when changing $\phi$ by $f$. 

The formula $(\ref{BWF})$ can be viewed as a natural extension of the Bochner-Weitzenb\"ock formula. In the literature, $Ric(L)=Ric+\nabla^2 \phi$ is called the infinite dimensional Bakry-Emery Ricci curvature associated with the Witten Laplacian $L$  on the weighted Riemannian manifold $(M, g, \phi)$. For $m\in [n, \infty]$,   we introduce 
\begin{eqnarray*}
Ric_{m, n}(L)=Ric+\nabla^2\phi-{\nabla\phi\otimes\nabla\phi\over m-n}, 
\end{eqnarray*}
and call it  the $m$-dimensional Bakry-Emery Ricci curvature associated with the Witten Laplacian $L$ on $(M, g, \phi)$.  Following Bakry and Emery \cite{BE}, we say that $(M, g, \phi)$ satisfies the $CD(K, m)$-condition if 
$$
Ric_{m, n}(L)\geq Kg.$$
Here we make the convention that $m=n$ if and only if $L=\Delta$,  $\phi$ is a constant, and  $Ric_{n, n}(\Delta)=Ric$.  Note that,  for $C^2$-smooth potential function $\phi$ on $(M, g)$  we have
$$Ric(L)=Ric_{\infty, n}(L)=\lim\limits_{m\rightarrow \infty}Ric_{m, n}(L).$$
Now  it is  well-known that the quantities $Ric(L)$ and $Ric_{m, n}(L)$  play as a good substitute of the Ricci curvature in many problems in comparison geometry and analysis on complete Riemannian manifolds with smooth weighted volume measures. See \cite{BE, BL,  Qian, Li05, Li12, Lich, Lot,  LoV, WW} and reference therein.

In the case of  Riemannian manifolds with a family of time dependent metrics and potentials, we call $(M, g(t), \phi(t), t\in [0, T])$ a $(K, m)$-super Perelman  Ricci flow if the metric $g(t)$ and the potential function $\phi(t)$ satisfy 
\begin{eqnarray}
{1\over 2}{\partial g\over \partial t}+Ric_{m, n}(L)\geq Kg, \label{KmsRF}
\end{eqnarray}
where $$L=\Delta_{g(t)}-\nabla_{g(t)}\phi(t)\cdot\nabla_{g(t)}$$ is the time dependent Witten Laplacian on $(M, g(t), \phi(t), t\in [0, T])$, and $K\in \mathbb{R}$ is a constant. When $m=\infty$, i.e., 
if the metric $g(t)$ and the potential function $\phi(t)$ satisfy the following inequality
\begin{eqnarray*}
{1\over 2}{\partial g\over \partial t}+Ric(L)\geq Kg,
\end{eqnarray*}
we call $(M, g(t), \phi(t), t\in [0, T])$   a $(K, \infty)$-super Perelman Ricci flow or a $K$-super Perelman Ricci flow. Indeed the $(K, \infty)$-Perelman Ricci flow (called also the $K$-Perelman Ricci flow)  
\begin{eqnarray*}
{1\over 2}{\partial g\over \partial t}+Ric(L)=Kg
\end{eqnarray*}
is a straightforward extension of the modified Ricci flow $\partial_t g=-2Ric(L)$ introduced by Perelman \cite{P1}  as the gradient flow of $\mathcal{F}(g, \phi)=\int_M (R+|\nabla \phi|^2)e^{-\phi}dv$ on $\mathcal{M}\times C^\infty(M)$ 
under the constraint condition that the measure $d\mu=e^{-\phi}dv$ is preserved. 

Super Ricci flows are super solutions to the Ricci flow. In \cite{MT}, McCann and Topping proved the equivalence between the  super Ricci flow 
$${\partial g \over \partial t}\geq -2Ric$$
and the contraction property of the $L^2$-Wasserstein distance between two solutions of the conjugate  heat equation
\begin{eqnarray*}
\partial_t u=\Delta u-{1\over 2}{\rm Tr}\left({\partial g \over \partial t}\right) u
\end{eqnarray*}
with different initial data. See also \cite{Lo2}.  In \cite{St1, St3},  K. T.  Sturm developed this idea to characterize  $(0, \infty)$-super Ricci flows on metric measure spaces.  In \cite{HN}, 
R. Haslhofer and A. Naber proved the  characterization of the super Ricci flows   $\partial_t g\geq -2Ric_g$ by  various functional inequalities and gradient estimates for  the heat equation $\partial_t u =\Delta_{g_t} u$ on $(M, g_t)$.

%
%

Since Perelman \cite{P1} introduced the $W$-entropy and proved its monotonicity for the Ricci flow, many
people have studied the $W$-like entropy for other geometric
flows on Riemannian manifolds \cite{N1, N2, LNVV, KN, LX, Li11, Li12}.  In \cite{N1, N2}, Ni proved an analogue of Perelman's $W$-entropy formula for the heat equation $\partial_t u=\Delta u$ on 
complete Riemannian manifolds with  fixed metric and with non-negative Ricci curvature.  In \cite{LX}, Li and Xu extended Ni's $W$-entropy formula to the heat equation $\partial_t u=\Delta u$ on complete Riemannian manifolds with fixed metric satisfying  $Ric\geq -Kg$, where $K\geq 0$ is a constant.  In \cite{Li11, Li12, Li16},  the second author of this paper proved the $W$-entropy formula  for  the Fokker-Planck equation and the heat equation of the Witten Laplacian on complete 
 Riemannian manifolds with the $CD(0, m)$-condition and gave a natural probabilistic interpretation of  Perelman's $W$-entropy for the Ricci flow.  In \cite{LL15}, we proved  the $W$-entropy formula for the heat equation of the Witten Laplacian on complete Riemannian manifolds with the  $CD(K, m)$-condition and on compact manifolds with  $(K, m)$-super Ricci flows, where $m\in [n, \infty)$ and $K\in \mathbb{R}$.  More precisely, let $(M, g(t), \phi(t), t\in [0, T])$ be a compact manifold with a $(K, m)$-super Ricci Perelman flow, $u={e^{-f}\over (4\pi t)^{m/2}}$  the fundamental solution to
the heat equation associated with the time dependent Witten Laplacian 
\begin{eqnarray}
\partial_t u=Lu. \label{Lu}
\end{eqnarray}
Define the $W$-entropy for the heat equation $(\ref{Lu})$ as follows
\begin{eqnarray}
W_{m, K}(u)={d\over dt}(tH_{m, K}(u))=H_{m, K}(u)+t{d\over dt}H_{m, K}(u), \label{WmK}
\end{eqnarray}
where
\begin{eqnarray}
H_{m, K}(u)=-\int_M u\log ud\mu-{m\over 2}\left(1+\log(4\pi t)+Kt+{K^2t^2\over 6}\right). \label{HmK}
\end{eqnarray}
In \cite{LL15}, we proved the following $W$-entropy formula 
\begin{eqnarray}
{d\over dt}W_{m, K}(u)&=&-2t\int_M \left[\left|\nabla^2f-\left({K\over 2}+{1\over 2t}\right)g\right|^2+\left({1\over 2}{\partial g\over \partial t}+Ric_{m, n}(L)+Kg\right)(\nabla f, \nabla f)\right] ud\mu\nonumber\\
& &\hskip1.5cm -{2t\over m-n}\int_M \left|\nabla \phi\cdot \nabla f+{(m-n)(1+Kt)\over 2t}\right|^2ud\mu. \label{WmKt}
\end{eqnarray}
In  time independent case, we pointed out in \cite{LL17b}  a close and deep connection between the Li-Yau-Hamilton type Harnack inequality and the $W$-entropy for the Witten Laplacian on complete Riemannian manifolds with the $CD(K, m)$-condition.

The purpose of this paper is to  establish the $W$-entropy  formula for the heat equation associated with the time dependent Witten Laplacian $L=\Delta_{g(t)}-\nabla_{g(t)}\phi(t)\cdot\nabla_{g(t)}$  on manifolds equipped 
 with $(K, \infty)$-super Perelman  Ricci flows. We would like to point out that we cannot use  the same definition formulas $(\ref{HmK})$ and  $(\ref{WmK})$ to  introduce the $W$-entropy on complete Riemannian manifolds with the $CD(K, m)$-condition and on manifolds with $(K, m)$-super Ricci flow. Indeed, when $m=\infty$ and $Ric(L)\geq Kg$ or $m=\infty$ and ${1\over 2}\partial_t g+Ric(L)\geq Kg$, neither the definition formula $(\ref{WmK})$  (resp.   $(\ref{HmK})$) for $H_{m, K}$ (resp. $W_{m, K}$)  nor the $W$-entropy formula $(\ref{WmKt})$ for $W_{m, K}$ make sense. 
  
To describe the idea how to introduce the $W$-entropy for the heat equation of the Witten Laplacian on manifolds with the $CD(K, \infty)$-condition and on manifolds with $(K, \infty)$-super Ricci flows, let us recall that  Bakry and Ledoux \cite{BL} proved the following characterization of Riemannian manifolds with the $CD(K, \infty)$-condition: Let  $(M, g)$  be a Riemannian manifold,  $\phi\in C^2(M)$, and $d\mu=e^{-\phi}dv$,  then the  $CD(K, \infty)$-condition holds,
 i.e., $Ric(L)\geq Kg$, if and only if the following reversal logarithmic Sobolev inequalities hold for the heat semigroup $P_t=e^{tL}$ generated by the Witten Laplacian $L=\Delta-\nabla\phi\cdot\nabla$ on $(M, g, \phi)$
\begin{eqnarray}
{|\nabla P_tf|^2\over P_tf}\leq {2K\over e^{2Kt}-1}(P_t(f\log f)-P_tf\log P_tf). \label{BLBL}
\end{eqnarray}
Inspired by $(\ref{BLBL})$, we introduce the revised Boltzmann-Shannon entropy $H_{K}(f, t)$ as follows
\begin{eqnarray*}
H_{K}(f, t)= D_K(t)({\rm Ent}(f|\mu)-{\rm Ent}(P_tf|\mu))
\end{eqnarray*}
where ${\rm Ent}(f|\mu)=\int_M f\log f d\mu$ is the  Boltzmann-Shannon entropy of the probability measure $f\mu$ with respect to the weighted volume measure $\mu$ on $(M, g)$, $D_0(t)={1\over t}$  and $D_{K}(t)={2K\over 1-e^{-2Kt}}$ for $K\neq 0$.  We then notice that if $Ric(L)\geq Kg$ then  Bakry and Ledoux's logarithmic Sobolev inequality yields
\begin{eqnarray*}
{d\over dt}H_{K}(f, t)\leq 0.
\end{eqnarray*}
Under the same condition $Ric(L)\geq Kg$, we prove that $H_{K}(f, t)$ satisfies a new second order differential inequality (see Theorem \ref{WCDK} below) 
\begin{eqnarray}
{d^2\over dt^2}H_K(f, t)+2K\coth(Kt) {d\over dt}H_K(f, t)+2D_K(t)\int_M |\nabla^2 \log P_tf|^2 P_tfd\mu\leq 0,\label{KKK}
\end{eqnarray}
and the equality in $(\ref{KKK})$ holds at some $t=t_0>0$ for non trivial $f$ if and only if $(M, g, \phi)$ is a gradient $K$-Ricci soliton, i.e., 
\begin{eqnarray}
Ric(L)=Kg.
\end{eqnarray}

 We now describe the main results of this paper. Our first result is the following theorem which extends above mentioned result due to  Bakry and Ledoux \cite{BL} to  manifolds with time dependent metrics and potentials.

\begin{theorem}\label{Thm-RLSI} Let $K\in\mathbb{R}$, $M$ be a manifold
equipped with a family of time dependent complete Riemannian metrics and
$C^2$-potentials $(g(t), \phi(t), t\in [0, T])$. Let
$L_t=\Delta_{g(t)}-\nabla_{g(t)}\phi(t)\cdot\nabla_{g(t)}$ be the time
dependent weighted Laplacian on $(M, g(t), \phi(t))$, $P_{s, t}$ be the time inhomogenuous heat semigroup generated by $L$, i.e., $u(s, t, x)=P_{s, t}f(x)$ is the solution to the heat equation $\partial_t
u=L u$ with the initial condition $u(s, s, \cdot)=f$, where $0\leq s<t\leq T$, and $f\in C(M,  (0, \infty))$. Then the following statements are equivalent:\\

(i) $(M, g(t), \phi(t), t\in [0,
T])$ is a $K$-super Perelman Ricci flow in the sense that
\begin{eqnarray}
\label{RFK} {1\over 2}{\partial g\over \partial t} +
Ric(L)\geq Kg, \label{Ksuper1}
\end{eqnarray}
 
(ii) for $0 \leq s < t \leq T$, the following logarithmic Sobolev inequality holds
\begin{eqnarray}
P_{s, t}(f\log f)-P_{s, t}f\log P_{s,t}f\leq {1-e^{-2K(t-s)}\over 2K}P_{s, t}\left({|\nabla f|^2\over
f}\right), \label{Ksuper2}
\end{eqnarray} 

(iii) for $0 \leq s < t \leq T$, the following reversal
logarithmic Sobolev inequality holds
\begin{eqnarray}
 {|\nabla P_{s,t} f|^2\over P_{s, t}
f}\leq {2K\over e^{2K(t-s)}-1}\left(P_{s, t}(f\log f)-P_{s,t}f\log P_{s, t}f\right),\label{Ksuper3}
\end{eqnarray}

(iv) for all $0\leq s<t\leq T$, the following Poincar\'e inequality holds
\begin{eqnarray}
P_{s,t}f^{2}- (P_{s,t}f)^{2} \leq \frac{1-e^{-2K(t-s)}} {K} P_{s, t}\left({|\nabla f|^2}\right), \label{Ksuper4}
\end{eqnarray}

(v) for all $0\leq s<t\leq T$, the following reversal  Poincar\'e inequality holds
\begin{eqnarray}
 {|\nabla P_{s,t} f|^2\over P_{s, t} f}\leq  {K\over e^{2K(t-s)}-1}(P_{s,t}f^{2}- (P_{s,t}f)^{2}). \label{Ksuper5}
\end{eqnarray}

(vi)  for $0 \leq s < t \leq T$, the following gradient estimate holds
\begin{eqnarray}
|\nabla P_{s, t}f|^2\leq e^{-2K(t-s)}P_{s, t}(|\nabla f|^2). \label{Ksuper6}
\end{eqnarray}

\end{theorem}

\begin{remark}{\rm  Theorem \ref{Thm-RLSI} can be viewed as a generalization of the well-known result due to Bakry and Ledoux \cite{BL} for the equivalence between the $CD(K, \infty)$-condition, the 
 logarithmic Sobolev inequalities, the Poincar\'e inequalities and the gradient estimate $(\ref{Ksuper6})$ for  the heat semigroup generated by the time independent Witten Laplacian on complete Riemannian manifolds. The proof of Theorem \ref{Thm-RLSI} is inspired by the semigroup argument due to Bakry and Ledoux \cite{BL}.  Indeed, using a similar approach as in the proof of Theorem \ref{Thm-RLSI}, we can further prove
the equivalence between the Bakry-Ledoux-Gromov-L\'evy isoperimetric 
inequality (see \cite{BL96})  and the $(K, \infty)$-super Perelman Ricci flows $(\ref{RFK})$. To save the length of the paper, we will do this in a forthcoming paper.   In \cite{St1, St3}, Sturm introduced the notion of $(0, \infty)$ and $(0, N)$-super Ricci flows on metric measure
spaces, and proved the equivalence between the $(0, \infty)$-super Ricci flows, the Poincar\'e inequality  and
the gradient estimate $(\ref{Ksuper6})$  for the heat semigroup $P_{s, t}$ generated by $L$ on metric measure
spaces. For $N<\infty$, the equivalence between the $(0, N)$-super Ricci flows and 
an improved version  the gradient estimate $(\ref{Ksuper6})$  for the heat semigroup $P_{s, t}$ generated by $L$ on metric measure
spaces was also proved in \cite{St3}.  When $\phi=0$ and $K=0$, R. Haslhofer and A. Naber  \cite{HN} proved the  characterization of the super Ricci flows   $\partial_t g\geq -2Ric_g$ by the Log-Sobolev inequality, the Poincar\'e inequality and the gradient estimate $(\ref{Ksuper6})$.  
We would like to mention that our work is independent of \cite{St1, St3, HN}. The first version of our work was posted on arxiv on  22 December 2014  and the second version was posted on arxiv  on 7 February 2016. See \cite{LL16}. 
 }
\end{remark}

As a byproduct of Theorem \ref{Thm-RLSI},  we derive the following Hamilton type dimension free Harnack inequality  for positive and bounded solution to the heat equation of the Witten Laplacian on complete Riemannian manifolds with $(K, \infty)$-super Perelman Ricci flows.  

\begin{theorem}\label{Thm0}\footnote{See also our paper \cite{LL17a} for a probabilistic proof of Theorem \ref{Thm0}.}
     Let $M$ be a manifold equipped with a family of time dependent complete Riemannian metrics and $C^2$-potentials
$(g(t), \phi(t), t\in [0, T])$  which is a $(-K, \infty)$-super Perelman
Ricci flow, i.e., 
\begin{eqnarray*}
{1\over 2}{\partial g\over \partial t}+Ric(L)\geq -Kg.
\end{eqnarray*}
where $K\geq 0$ is a constant independent of  $t\in [0, T]$. Let $u$ be a positive and bounded solution to the heat equation
\begin{eqnarray*}
\partial_t u=Lu,
\end{eqnarray*}
where
$$L=\Delta_{g(t)}-\nabla_{g(t)}\phi(t)\cdot\nabla_{g(t)}
$$
is the time dependent Witten Laplacian on $(M, g(t), \phi(t))$.
 Then the sharp version of the Hamilton Harnack inequality holds: for all $x\in M$ and $t>0$,
\begin{eqnarray}
{|\nabla u|^2\over u^2}\leq {2K\over 1-e^{-2Kt}}\log (A/u),\ \ \label{HH}
\end{eqnarray}
where
\begin{eqnarray*}
A:=\sup\limits\{u(t, x): x\in M, t\geq 0\}.
\end{eqnarray*}
In particular, the Hamilton Harnack inequality holds 
\begin{eqnarray}
{|\nabla u|^2\over u^2}\leq \left({1\over
t}+2K\right)\log(A/u).\label{Ham}
\end{eqnarray}
\end{theorem}
%

Now we introduce the $W$-entropy for the heat equation of the Wtiietn Lalpacian on manifolds with the $CD(K, \infty)$-condition as follows
\begin{eqnarray}
W_K(f, t)=H_K(f, t)+{\sinh(2Kt)\over 2K}{d\over dt}H_K(f, t).\label{Wk}
\end{eqnarray} 
The definition formula $(\ref{Wk})$ is new and is different from Perelman's $W$-entropy $(\ref{entropy-1})$ for the Ricci flow and the $W$-entropy $(\ref{WmK})$ for  the heat equation of the Witten Lalpacian on manifolds with the $CD(K, \infty)$-condition.   Nevertheless, we have the following $W$-entropy formula for 
the heat equation $\partial_t u=Lu$ on manifolds with 
the $CD(K, \infty)$-condition and  on $K$-super Perelman Ricci flows. Our result gives  a new characterization of  the $K$-Ricci soliton  and the $(K, \infty)$-Perelman Ricci flow.

\begin{theorem} \label{WCDK1} Let $(M, g, \phi)$ be a complete Riemannian manifold with bounded geometry condition and $Ric(L)\geq Kg$, where $K\in \mathbb{R}$ is
a constant.  Then 
\begin{eqnarray*}
{d\over dt}W_{K}(f, t)=-(1+e^{2Kt})\int_M \left(|\nabla^2 \log P_tf|^2+(Ric(L)-Kg)(\nabla\log P_tf, \nabla\log P_tf)\right)P_tfd\mu.\label{WKW1}
\end{eqnarray*}
In particular, if $Ric(L)\geq Kg$, we have 
\begin{eqnarray*}
{d\over dt}W_{K}(f, t)+(1+e^{2Kt})\int_M |\nabla^2 \log P_tf|^2P_tfd\mu\leq 0,
\end{eqnarray*}
and the equality hods at some time $t=t_0>0$ if and only if $(M, g, \phi)$ is a gradient $K$-Ricci soliton, i.e.,  
 \begin{eqnarray*}
 Ric+\nabla^2\phi=Kg.
 \end{eqnarray*}
In particular, this is the case when $L=\Delta-Kx\cdot\nabla$ is  the Ornstein-Ulenbeck operator on the Gaussian space over $\mathbb{R}^n$. 

\end{theorem}

\begin{theorem} \label{WCDK2} Let $(M, g(t), \phi(t), t\in [0, T])$ be a compact manifold with a family of time dependent metrics and potentials such that
 \begin{eqnarray}
{1\over 2}{\partial g\over \partial t}+Ric(L)\geq Kg, \ \ \ {\partial\phi\over \partial t}={1\over 2}{\rm Tr}\left( {\partial g\over \partial t}\right). \label{HHH9}
\end{eqnarray}
where  $K\in \mathbb{R}$. Then, for all $t\in [0, T]$, we have
\begin{eqnarray*}
{d\over dt}W_{K}(f, t)=-(1+e^{2Kt})\int_M \left({1\over 2}\left[{\partial g\over \partial t}+\Gamma_2\right]-Kg\right) (\nabla\log P_tf, \nabla\log P_tf)P_tfd\mu.\label{WWWW}
\end{eqnarray*}
where
 $$\Gamma_2(\nabla u, \nabla u)=2|\nabla^2 u|^2+2Ric(L)(\nabla u, \nabla u).$$ In particular, we have
\begin{eqnarray*}
{d\over dt}W_{K}(f, t)+(1+e^{2Kt})\int_M |\nabla^2\log P_tf|^2 P_t fd\mu\leq 0, 
\end{eqnarray*}
 and the equality holds on $(0, T]$  if and only if $(M, g(t), \phi(t), t\in [0, T])$ is the $K$-Perelman Ricci flow 
\begin{eqnarray*}
{1\over 2}{\partial g\over \partial t}+Ric(L)=Kg,  \ \ \ {\partial\phi\over \partial t}=-R-\Delta \phi+nK.
\end{eqnarray*}
\end{theorem}

In particular, for $K=0$, we have the following 
\begin{corollary}  Let $(M, g(t), \phi(t), t\in [0, T])$ be a compact manifold with a family of time dependent metrics and potentials such that 
\begin{eqnarray*}
{1\over 2}{\partial g\over \partial t}+Ric(L)\geq 0,  \ \ \ {\partial\phi\over \partial t}={1\over 2}{\rm Tr}\left( {\partial g\over \partial t}\right).
\end{eqnarray*}
Then
\begin{eqnarray*}
{d\over dt}W_{0}(f, t)=-2\int_M \left({1\over 2}\left[{\partial g\over \partial t}+\Gamma_2\right]\right) (\nabla\log P_tf, \nabla\log P_tf)P_tfd\mu.\label{WWWW}
\end{eqnarray*}
In particular, for all $t\in (0, T]$, we have 
\begin{eqnarray*}
{d\over dt}W_{0}(f, t)+2\int_M |\nabla^2\log P_tf|^2 P_t fd\mu\leq 0,
\end{eqnarray*}
and the equality holds on $(0, T]$ if and only if $(M, g(t), \phi(t), t\in [0, T])$ is the Perelman Ricci flow
\begin{eqnarray*}
{\partial g\over \partial t}=-2(Ric+\nabla^2\phi),  \ \ \ {\partial\phi\over \partial t}=-R-\Delta \phi,
\end{eqnarray*}
which is the gradient flow of $\mathcal{F}(g, \phi)=\int_M[R+|\nabla \phi|^2]e^{-\phi}dv$ under the constraint the measure $d\mu=e^{-\phi}dv$ is preserved. 
\end{corollary}

Theorem \ref{WCDK1} and  Theorem \ref{WCDK2} have been announced in our survey  paper \cite{LL17c}.

The rest of this paper is organized as follows.  In Section $2$, we prove Theorem \ref{Thm-RLSI} and Theorem \ref{Thm0}.   In Section $3$,  we prove Theorem \ref{WCDK1}. In Section $4$, we prove Theorem \ref{WCDK2}. 
In Section $5$, we prove a second order differential entropy inequality on $(K, m)$-super Ricci flows, which can be used to characterize the $(K, m)$-Ricci solitons and the $(K, m)$-Ricci flows.  Finally, we give a probabilistic interpretation of the $W$-entropy for the heat equation of the Witten Laplacian on manifolds with the $CD(K, m)$-condition.

To end this section, let us mention that this paper is a revised version of  a part of our 2014-2016 preprint \cite{LL16},  which contained also the Li-Yau and the Li-Yau-Hamilton type Harnack inequalities on variants of $(K, m)$-super Perelman  Ricci flows. As the preprint \cite{LL16} is too long, we  have divided it into three papers (see also \cite{LL17a, LL17b}).  

\section{Log-Sobolev inequalities and Harnack inequalities for Witten Laplacian}

\subsection{Log-Sobolev inequalities and $K$-super Perelman  Ricci flows}

\noindent{\bf Proof of Theorem \ref{Thm-RLSI}}.  Note that $\partial_t P_{s, t} f=L_t P_{s, t}f$, and $\partial_s P_{s, t}f=-P_{s, t}L_s f$.
Let 
\begin{eqnarray*}
h(s, t)=e^{-2Kt}P_{s+T-t, T}\left({|\nabla P_{s, s+T-t}f|^2\over P_{s, s+T-t}f}\right),
\ \ \ t\in [s, T].
\end{eqnarray*}
Note that,  at time $T-t+s$, the generalized Bochner formula implies
\begin{eqnarray}
\label{Boh1}
(\partial_t+L){|\nabla u|^2\over u}={2\over u}|\nabla^2
u-u^{-1}\nabla u\otimes \nabla u|^2+2u^{-1}\left({1\over 2}{\partial
g\over
\partial t}+Ric(L)\right)(\nabla u, \nabla u).
\end{eqnarray}
Differentiating $h(s, t)$ with respect to $t$ in $[s, T]$, we have
\begin{eqnarray}
& &\partial_t h(s, t)=-2K h(s, t)+e^{-2Kt}P_{s+T-t, T}\left[\left({\partial \over
\partial t}+L\right)\left({|\nabla P_{s, s+T-t}f|^2\over
P_{s, s+T-t}f}\right)\right]\nonumber \\
&=&e^{-2Kt} P_{s+T-t, T}\left[{2\over u}\left|\nabla^2 u-{\nabla u\otimes \nabla u\over u}\right|^2+{2\over u}\left({1\over 2}{\partial g\over \partial t}+Ric(L)-Kg\right)(\nabla u, \nabla u) \right]\nonumber\\
&\geq& 2e^{-2Kt} P_{s+T-t, T}\left[u^{-1}\left({1\over 2}{\partial g\over
\partial t}+Ric(L)-Kg\right)(\nabla u, \nabla u)\right].\label{hst1}
\end{eqnarray}

Assuming $(i)$ holds, i.e.,  $(M, g(t), \phi(t))$ is a $(K, \infty)$-super Perelman Ricci flow, we have
\begin{eqnarray}
\partial_t h(s, t) \geq 0.\label{hst2}
\end{eqnarray}
Thus, $t\rightarrow h(s, t)$ is increasing on $[s, T]$. This yields,
for all $t\in (s, T)$, 
\begin{eqnarray}
e^{2K(t-s)}{|\nabla P_{s, T}f|^2\over P_{s, T}f}\leq P_{s+T-t, T}\left({|\nabla
P_{s, s+T-t}f|^2\over P_{s, s+T-t}f}\right)\leq e^{-2K(T-t)}P_{s,T}\left({|\nabla
f|^2\over f}\right). \label{hst3}
\end{eqnarray}
Fix $s$ and $T$, differentiating $\alpha(r):=P_{r, T}(P_{s, r}f\log
P_{s, r}f)$ with respect to $r$ in $[s, T]$, using $\partial_r P_{s, r}f=L_r P_{s, r}f=-P_{s, r}L_s f$ and the chain rule, we have
\begin{eqnarray*}
\alpha'(r)&=&P_{r, T}\left((-L_{r}+\partial_r)(P_{s, r}f\log
P_{s, r}f)\right)\\
&=&P_{r, T}\left[-L_rP_{s, r}f \log P_{s, r}f-{2|\nabla P_{s, r}f|^2\over P_{s, r}f}-P_{s, r}f \left({LP_{s, r}f\over P_{s, r}f}-{|\nabla P_{s, rf}|^2\over P_{s, r}f}\right)\right]\\
& & \hskip2cm +P_{r, T}\left(L_r P_{s, r}f \log P_{s, r}f+\partial_r P_{s, r}f\right)\\
&=& -P_{r, T}\left({|\nabla P_{s, r} f|^2\over P_{s, r}f}\right).
\end{eqnarray*}
Hence
\begin{eqnarray*}
{d\over dt}P_{s+T-t, T}(P_{s, s+T-t}f\log
P_{s, s+T-t}f)
=P_{s+T-t, T}\left({|\nabla P_{s, s+T-t} f|^2\over P_{s, s+T-t}f}\right).
\end{eqnarray*}
Integrating in $t$ from $s$ to $T$ and using $(\ref{hst3})$, we have
\begin{eqnarray*}
& &P_{s,T}(f\log f)-P_{s,T}f\log P_{s, T} f
=\int_s^T P_{s+T-t, T}\left({|\nabla P_{s, s+T-t} f|^2\over P_{s, s+T-t}f}\right)dt\\
& &\hskip1cm  \leq \int_s^T e^{-2K(T-t)} P_{s,T}\left({|\nabla
f|^2\over f}\right)dt= \frac{1-e^{-2K(T-s)}} {2K}P_{s, T}\left({|\nabla f|^2\over f}\right).
\end{eqnarray*}
Similarly, we have
\begin{eqnarray*}
& &P_{s, T}(f\log f)-P_{s, T}f\log P_{s,T} f=\int_s^T P_{s+T-t, T}\left({|\nabla P_{s, s+T-t} f|^2\over P_{s, s+T-t}f}\right)dt\nonumber\\
& &\hskip1cm \geq \int_s^T e^{2K(t-s)}{|\nabla P_{s,T} f|^2\over
P_{s,T} f}dt={e^{2K(T-s)}-1\over 2K}{|\nabla P_{s,T} f|^2\over P_{s, T} f}.\nonumber
\end{eqnarray*}
Changing $T$ by $t$, we obtain $(\ref{Ksuper2})$ and $(\ref{Ksuper3})$. This proves that $(i)$ implies $(ii)$ and $(iii)$. 

Assuming $(iii)$ holds, applying the log-Sobolev
inequality $(\ref{Ksuper2})$  to
$1+\varepsilon f$ and using the Taylor expansion $\log(1+\varepsilon f)=\varepsilon f-{\varepsilon^2f^2\over 2}+o(\varepsilon^2)$ for $\varepsilon\rightarrow 0$, we have
\begin{eqnarray*}
{\varepsilon^2\over 2}(P_{s, t}f^2-(P_{s, t}f)^2)+o(\varepsilon^2)\leq \frac{1-e^{-2K(T-s)}} {2K} \varepsilon^2 P_{s, t}(|\nabla f|^2)+o(\varepsilon^2).
\end{eqnarray*}
This yields the Poincar\'e inequality $(\ref{Ksuper4})$, i.e., 
\begin{eqnarray*}
P_{s,t}f^{2}- (P_{s,t}f)^{2} \leq \frac{1-e^{-2K(t-s)}} {K} P_{s, t}\left({|\nabla f|^2}\right).
\end{eqnarray*}
So $(ii)$ implies $(iv)$. Similarly, applying the reversal log-Sobolev
inequality  $(\ref{Ksuper3})$ to
$1+\varepsilon f$ and using Talyor expansion for $\varepsilon\rightarrow 0$, we
obtain the reversal Poincar\'e inequality $(\ref{Ksuper5})$, i.e., 
\begin{eqnarray*}
{|\nabla P_{s,t} f|^2\over P_{s, t} f}\leq  {K\over e^{2K(t-s)}-1}(P_{s,t}f^{2}- (P_{s,t}f)^{2}). 
\end{eqnarray*}
This proves that $(iii)$ implies $(v)$. 

To prove that $(iv)$ implies $(i)$, set
\begin{eqnarray*}
w(s, t)=P_{s,t}f^{2}- (P_{s,t}f)^{2}+\frac{e^{2K(s-t)}-1}{K} P_{s, t}\left({|\nabla f|^2}\right).
\end{eqnarray*}
Taking derivatives in $s$ , we have
\begin{eqnarray*}
\partial_{s} w(s, t) &=& -P_{s, t}L_sf^{2} + 2P_{s, t}f P_{s, t}L_sf +2e^{2K(s-t)}P_{s, t} (|\nabla f|^{2}) \\
& &\ \ \ + \frac{e^{2K(s-t)} -1} {K} [P_{s, t}(-L_s|\nabla f|^{2}+\partial_s|\nabla f|^{2}_{g(s)})],
\end{eqnarray*}
and
\begin{eqnarray*}
\left.\partial_{s}^2 w(s, t)\right|_{s=t} 
&=&[-\partial_s(P_{s, t}L_sf^{2}) + 2\partial_s(P_{s, t}fP_{s, t}L_sf) +4Ke^{2K(s-t)}P_{s, t}(|\nabla f|^{2})]_{s=t} \\
& &\ \ + [2e^{2K(s-t)}\partial_sP_{s, t}(|\nabla f|^{2}) + 2e^{2K(s-t)}P_{s, t}(-L_s(|\nabla f|^{2}) + \partial_s|\nabla f|^{2}_{g(s)})]_{s=t}\\
&=& L^{2}_sf^{2}  - \partial_s L_sf^{2}- 2(L_s f)^{2} - 2fL^{2}_{s}f + 2f\partial_sL_sf\\\
& &\ \ \ \ + 4K |\nabla f|^{2} - 4L_s|\nabla f|^{2}
+4\partial_s|\nabla f|_{g(s)}^2.
\end{eqnarray*}
Note that  $w(t, t)=0$ and $\left. \partial_s w(s,
t)\right|_{s=t}=0$. As for all $s<t$, $w(s, t)\leq 0$, the Taylor expansion yields $\left.\partial^2_s w(s, t)\right|_{s=t}\leq 0$. 
Using the fact
\begin{eqnarray*}
\partial_s |\nabla f|_{g(s)}^2&=&-\partial_s g(s)(\nabla f, \nabla f),\\
\partial_s L_s f^2&=&\partial_s (2fL_s f+2|\nabla f|_{g(s)}^2)=2f\partial_s L_s f+2\partial_s |\nabla f|_{g(s)}^2,\\
L_s^2 f^2&=&L_s(2fL_sf+2|\nabla f|_{g(s)}^2)\\
&=&2(L_s f)^2+2fL_s^2f+4\nabla f\cdot\nabla L_sf+2L_s|\nabla
f|_{g(s)}^2,
\end{eqnarray*}
and by the  generalized Bochner formula, we have 
\begin{eqnarray*}
\left. \partial_{s}^2 w(s, t)\right|_{s=t}
&=&-4\left[\|\nabla^2 f|_{g(t)}^2+Ric(L_t)(\nabla f, \nabla f)\right]+4K|\nabla f|^2_{g(t)}-2\partial_t g(t)(\nabla f, \nabla f)\\
&=&-4\left[\|\nabla^2 f|_{g(t)}^2+\left({1\over 2}\partial_t
g(t)+Ric(L_t)-Kg\right)(\nabla f, \nabla f)\right].
\end{eqnarray*}
Taking $f$ to be normal coordinate functions near any fixed point $x$ on $(M,
g(t))$, we derive that, at any time $t \in [0, T]$,
$$
{1\over 2}{\partial g\over
\partial t} + Ric(L_t)\geq  Kg.
$$
So $(iv)$ implies $(i)$. Similarly, we can prove that $(iv)$ implies  $(i)$. 

Set $\psi(t)=e^{-2K(T-t))} P_{t, T}(|\nabla P_{s, t}f|^2$. Differentiating $\psi(t)$, we have
\begin{eqnarray*}
\psi'(t)=2K e^{-2K(T-t))}P_{t, T}(|\nabla P_{s, t}f|^2)+e^{-2K(T-t))}P_{t, T}(\partial_t-L_t)|\nabla P_{s, t}f|^2).
\end{eqnarray*} 
By the generalized Bochner formula, we have
\begin{eqnarray*}
(\partial_t-L_t) |\nabla P_{s, t}f|^2 =-2 |\nabla^2 P_{s, t}f|^2-2 \left( {1\over 2} {\partial g\over \partial t}+ Ric(L)\right) (\nabla P_{s, t}f, \nabla P_{s, t}f).
\end{eqnarray*}  
Thus
\begin{eqnarray*}
\psi'(t)=-2|\nabla^2 P_{s, t}f|^2-2\left({1\over 2}{\partial g\over \partial t}+Ric(L)-Kg\right)(\nabla P_{s, t}f, \nabla P_{s, t}f).
\end{eqnarray*} 
It follows from  $(\ref{Ksuper1})$ that $\psi'(t)\leq 0$ and $\psi(t)$ is decreasing on $[s, T]$, which yields $(\ref{Ksuper6})$, i.e., 
\begin{eqnarray*}
|\nabla P_{s, T}f|^2\leq e^{-2K(T-s)}P_{s, T}(|\nabla f|^2).
\end{eqnarray*}
Conversely, if $(\ref{Ksuper6})$ holds, we have $\psi'(t)\leq 0$ on $[s, T]$. In particular, $\psi'(s)\leq 0$. Thus
\begin{eqnarray*}
|\nabla^2 f|^2+\left({1\over 2}{\partial g\over \partial t}+Ric(L)-Kg\right)(\nabla f, \nabla f)\geq 0.
\end{eqnarray*} 
which yields $(\ref{Ksuper1})$ by taking special $f$ as the normal coordinate function at any fixed $x\in M$.  
The proof of Theorem \ref{Thm-RLSI}  is completed.   \hfill $\square$

\subsection{Hamilton's Harnack inequality on complete super Perelman Ricci flows}

In this subsection, we prove Hamilton's Harnack inequality for the time dependent Witten Laplacian on complete $(-K, \infty)$-super Perelman Ricci flows.

\medskip

\noindent{\bf Proof of Theorem \ref{Thm0}}. We modify  the method used in \cite{Li16}. Let $t\in [0, T)$  and $s\in [0, T-t]$. Using the reversal logarithmic Sobolev inequality and the fact $0<f\leq A$, we have
\begin{eqnarray*}
{|\nabla P_{s, s+t} f|^2\over P_{s, s+t} f}&\leq& {2K\over
1-e^{-2Kt}}\left(P_{s, s+t}(f\log f)-P_{s, s+t}f \log P_{s, s+t}f\right)\\
&\leq& {2K\over
1-e^{-2Kt}}\left(P_{s, s+t}(f\log A)-P_{s, s+t}f \log P_{s, s+t}f\right).
\end{eqnarray*}
Thus
\begin{eqnarray*}
|\nabla \log P_{s, s+t} f|^2\leq {2K\over 1-e^{-2Kt}}\log (A/P_{s, s+t}f).
\end{eqnarray*}
Using ${1\over 1-e^{-x}}\leq 1+{1\over x}$ for $x\geq 0$, we have
\begin{eqnarray*}
|\nabla \log P_{s, s+t} f|^2 \leq \left(2K+{1\over t}\right)\log
(A/P_{s, s+t}f).
\end{eqnarray*}
In particular, for $s=0$, we have 
\begin{eqnarray*}
{|\nabla u|^2\over u^2} \leq \left(2K+{1\over t}\right)\log
(A/u).
\end{eqnarray*}
The proof of Theorem \ref{Thm0} is completed. \hfill $\square$

\section{$W$-entropy for Witten Laplacian on manifolds with $(K, \infty)$-condition}

\subsection{A new second order differential inequality on Boltzmann-Shannon entropy}

In this subsection,  we  prove a new second order differential inequality for the  Boltzmann-Shannon entropy on complete Riemannian manifolds 
with  the $CD(K, \infty)$ condition.

Let $C_0(t)={1\over t}$ and $C_K(t)={2K\over e^{2Kt}-1}$  for $K\neq 0$.  Let  $D_0(t)={1\over t}$ and $D_{K}(t)={2K\over 1-e^{-2Kt}}$  for $K\neq 0$.
Then $D_K'(t)=-C_K(t)D_K(t)$ for all $K\in \mathbb{R}$ and $t>0$.
We first introduce the revised relative Boltzmann-Shannon entropy
\begin{eqnarray*}
H_{K}(f, t)=D_K(t)\int_M (f\log f-P_t f\log P_tf)d\mu,
\end{eqnarray*}
where $f\in C(M, (0, \infty))$.

\begin{theorem}\label{WCDK}  Let $(M, g)$ be a complete Riemannian manifold with bounded geometry condition, and $\phi\in C^4(M)$ with $\nabla\phi\in C^3_b(M)$. Suppose that $Ric(L)\geq Kg$, where $K\in \mathbb{R}$ is
a constant. Then
\begin{eqnarray*}
{d\over dt}H_{K}(f, t)\leq 0,\ \ \ \forall t>0,
\end{eqnarray*}
and for all $t>0$, we have
\begin{eqnarray}
& &{d^2\over dt^2}H_{K}(f, t)+2K\coth(Kt){d\over dt}H_K(t)+2D_K(t) \int_M  |\nabla^2
\log P_tf|^2 P_tf d\mu\nonumber\\
& &\hskip1.5cm =-2D_K(t)\int_M (Ric(L)-Kg)(\nabla \log P_tf, \nabla \log P_tf)P_tfd\mu. \label{HHH2}
\end{eqnarray}
In particular,  for all $t>0$, we have
\begin{eqnarray}
{d^2\over dt^2}H_K(f, t)+2K\coth(Kt) {d\over dt}H_K(f, t)+2D_K(t)\int_M |\nabla^2 \log P_tf|^2 P_tf d\mu\leq 0,\label{KK2}
\end{eqnarray}
and the equality in $(\ref{KK2})$ holds at some $t=t_0>0$ if and only if $Ric(L)=Kg$,  i.e., $(M, g, \phi)$ is a gradient $K$-Ricci soliton
\begin{eqnarray*}
Ric+\nabla^2\phi=Kg.
\end{eqnarray*}
\end{theorem}
{\it Proof}. By direct calculation and using the integration
by parts formula, we have
\begin{eqnarray}
{d\over d t} H_{K}(f, t)&=&-C_K(t) H_K(f, t)+D_K(t)\int_M {|\nabla P_t f|^2\over P_t f} d\mu\nonumber\\
&=&D_K(t)\int_M \left[{|\nabla P_tf|^2\over P_tf}+C_K(t)(P_tf\log P_tf-P_t(f\log f))\right]d\mu.\label{HHHH1}
\end{eqnarray}
Under the condition  $Ric(L)\geq Kg$, by the reversal  logarithmic Sobolev inequality due to Bakry and Ledoux \cite{BL}, for all $t>0$, we have
\begin{eqnarray*}
{|\nabla P_tf|^2\over P_tf}\leqq C_K(t)(P_t(f\log f)-P_tf\log P_tf). \label{BLRLSI}
\end{eqnarray*}
Hence, for all $K\in \mathbb{R}$, we have
\begin{eqnarray*}
{d\over d t} H_{K}(f, t)\leq 0,\ \ \ \forall\ t>0.
\end{eqnarray*}
Taking the time derivative on the both sides of $(\ref{HHHH1})$, we have
\begin{eqnarray*}
{d^2\over dt^2}H_{K}(f, t)
&=& D_K(t){d\over dt}\int_M {|\nabla
P_tf|^2\over P_tf}d\mu-C_K(t){d\over dt} H_K(f, t)\\
& &\ \ +{d\over dt}D_K(t)\int_M {|\nabla P_tf|^2\over P_tf}d\mu-{d\over dt} C_K(t)H_K(f, t).
\end{eqnarray*}
Let $u=P_tf$. By  \cite{Li12, Li16}, we have
\begin{eqnarray*}
{d\over dt}\int_M {|\nabla u|^2\over u}d\mu=-2\int_M \Gamma_2(\nabla \log u, \nabla \log u)ud\mu,
\end{eqnarray*}
where 
\begin{eqnarray*}
\Gamma_2(\nabla \log u, \nabla \log u)= |\nabla^2 \log u|^2+Ric(L)(\nabla \log u, \nabla \log u ).
\end{eqnarray*}
Note that
\begin{eqnarray*}
{d\over dt}C_K(t)={d\over dt}{2K\over e^{2Kt}-1}=-{2K\over 1-e^{-2Kt}}C_k(t)=-D_K(t)C_K(t),
\end{eqnarray*}
and
\begin{eqnarray*}
2K+C_K(t)=2K+{2K\over e^{2Kt}-1}={2K\over 1-e^{-2Kt}}=D_K(t).
\end{eqnarray*}
Hence
\begin{eqnarray}
{d^2\over dt^2}H_{K}(f, t)
&=&-2D_K(t) \int_M [\Gamma_2(\nabla \log u, \nabla \log u)-K|\nabla\log u|^2]ud\mu-C_K(t){d\over dt}H_K(f, t)\nonumber\\
& & -[2K+C_K(t)] D_K(t)\int_M |\nabla \log  u|^2 ud\mu+C_K(t)D_K(t)H_K(f, t)\nonumber\\
&=&-2D_K(t) \int_M (\Gamma_2-Kg)(\nabla \log u, \nabla\log u)ud\mu-C_K(t){d\over dt}H_K(f, t)\nonumber\\
& & -D^2_K(t)\int_M |\nabla \log  u|^2 ud\mu+C_K(t)D_K(t)H_K(f, t)\label{KH1}.
\end{eqnarray}
Combining $(\ref{HHHH1})$  with $(\ref{KH1})$,  we then finish the proof of Theorem \ref{WCDK}. \hfill $\square$ 

\subsection{$W$-entropy   for Witten Laplacian with $CD(K, \infty)$ condition }

In this subsection we introduce the $W$-entropy for the Witten Laplacian on manifolds satisfying the $CD(K, \infty)$ condition and  prove Theorem \ref{WCDK1}.  For this purpose,  let $\alpha_K: (0, \infty)\rightarrow (0, \infty)$ be a $C^1$-smooth function which will be determined later. Define the $W$-entropy by the revised Boltzmann entropy formula
\begin{eqnarray*}
W_{K}(f, t):={1\over \dot \alpha_K(t)}{d\over dt}(\alpha_K(t)H_K(f, t))= H_K+{\alpha_K\over \dot \alpha_k} \dot H_K.
\end{eqnarray*}
Set $\beta_K={\alpha_K\over \dot\alpha_K}$. Then
\begin{equation*}
{d\over dt}W_{K}(f, t)=\beta_K(\ddot{H}_K+\frac{1+\dot{\beta}_K}{\beta_K}\dot H_K).
\end{equation*}
Solving the ODE
\begin{equation*}
\frac{1+\dot{\beta_K}}{\beta_K}=2K\coth(Kt),
\end{equation*}
we have a special solution 
\begin{eqnarray*}
\beta_K(t)={\sinh(2Kt)\over 2K},
\end{eqnarray*}
and 
\begin{eqnarray*}
\alpha_K(t)=K\tanh(Kt).
\end{eqnarray*}
Therefore
\begin{eqnarray*}
W_K(f, t)=H_K(f, t)+{\sinh(2Kt)\over 2K}{d\over dt}H_K(f, t).
\end{eqnarray*}

\medskip

\noindent{\bf Proof of Theorem \ref{WCDK1}}. 
By Theorem \ref{WCDK}, we have
\begin{eqnarray*}
{d\over dt}W_{K}(f, t)&=&-{\sinh(2Kt)\over K}D_K(t)\int_M  |\nabla^2
\log P_tf|^2 P_tf d\mu\nonumber\\
& & -{\sinh(2Kt)\over K}D_K(t)\int_M (Ric(L)-Kg)(\nabla\log P_tf, \nabla\log P_tf)P_tfd\mu.\label{WWWW}
\end{eqnarray*}
Note that 
\begin{eqnarray*}
{\sinh(2Kt)\over K}D_K(t)={\sinh(2Kt)\over K} {2K\over 1-e^{-2Kt}}=1+e^{2Kt}.
\end{eqnarray*}
Thus
\begin{eqnarray*}
{d\over dt}W_{K}(f, t)=-(1+e^{2Kt})\int_M  [|\nabla^2
\log P_tf|^2 +(Ric(L)-Kg)(\nabla\log P_tf, \nabla\log P_tf)]P_tfd\mu.\label{WWWW}
\end{eqnarray*}

In particular, when $Ric(L)\geq  Kg$, then for all $t>0$, we have
\begin{eqnarray*}
{d\over dt}W_{K}(f, t)\leq -(1+e^{2Kt})\int_M  |\nabla^2 \log P_tf|^2 P_tf d\mu.
\end{eqnarray*}
Moreover, under the condition $Ric(L)\geq Kg$, we see that
\begin{eqnarray*}
{d\over dt}W_{K}(f, t)+(1+e^{2Kt})\int_M  |\nabla^2
\log P_tf|^2P_tfd\mu=0 \label{OUsoliton1}
\end{eqnarray*}
holds for non trivial $f$  at some $t=t_0>0$ if and only if $(M, g, \phi)$ is a gradient $(K, \infty)$-Ricci soliton
\begin{eqnarray*}
Ric(L)=Kg.\label{WWWW}
\end{eqnarray*}
This finishes the proof of Theorem \ref{WCDK1}.  \hfill $\square$

\subsection{Rigidity model for the $W_K$-entropy: gradient Ricci solitons}
 
 By Theorem \ref{WCDK1},  the rigidity model for the $W_K$-entropy on compact or complete Riemannian manifolds with $CD(K, \infty)$-condition  is  the gradient $(K, \infty)$-Ricci soliton
 $$Ric(L)=Kg,$$ 
 equivalently
 $$Ric+\nabla^2 \phi=Kg.$$ 
 
 By Hamilton \cite{Ha} and Ivey \cite{Iv1},  steady (i.e., $K=0$) or expanding (i.e., $K<0$) compact Ricci solitons must be trivial. That is to say, if $(M, g, \phi)$ is a compact Riemannian manifold with $Ric(L)=Ric+\nabla^2\phi=Kg$, where $K\leq 0$, 
 then $(M, g)$ must be Einstein with $Ric=Kg$ and $\phi$ must be constant.   
 
 The Gaussian soliton, i.e., $M=\mathbb{R}^n$ with Euclidean metric $g_0$,  and $\phi_K(x)={K\|x\|^2\over 2}+{n\over 2}\log(2\pi K^{-1})$, where $K>0$,  is a complete Ricci soliton with $Ric(L)=Kg_0$.   In this case, $d\mu(x)={1\over (2\pi K^{-1})^{n/2}} e^{-{K\|x\|^2\over 2}}dx$ is the Gaussian measure on $\mathbb{R}^n$, $L=\Delta-Kx\cdot \nabla$ is the Ornstein-Uhlenbeck operator on $\mathbb{R}^n$, and $Ric(L)=\nabla^2\phi_K=Kg_0$. See \cite{Iv2} for more examples of complete Ricci solitons.

\section{$W$-entropy for Witten Laplacian on compact $K$-super Perelman  Ricci flows}

In this section, we prove the $W$-entropy formula for the heat equation of the time dependent
Witten Laplacian on compact manifold equipped with a $K$-super Perelman  Ricci flow.

Let $(M, g(t), \phi(t), t\in [0, T])$ be a compact Riemannian manifold with a family of time dependent metrics $g(t)$ and potentials $\phi(t)$. Let
$$L=\Delta_{g(t)}-\nabla_{g(t)}\phi(t)\cdot\nabla_{g(t)}$$
be the time dependent Witten Laplacian on $(M, g(t), \phi(t))$. Let
$$d\mu(t)=e^{-\phi(t)}dvol_{g(t)}.$$
Suppose that
\begin{eqnarray}
{\partial \phi\over \partial t}={1\over 2}{\rm Tr} {\partial g\over
\partial t}.\label{PPPPP}
\end{eqnarray}
Then $\mu(t)$ is indeed independent of $t\in [0, T]$, i.e.,
\begin{eqnarray*}
{\partial \mu(t)\over \partial t}=0, \ \ \ t\in [0, T].
\end{eqnarray*}

\begin{theorem}\label{WCDK3} Let $(M, g(t), \phi(t), t\in [0, T])$ be a compact $K$-super Perelman  Ricci flow, i.e., 
$${1\over 2}{\partial g\over \partial t}+Ric(L)\geq Kg,$$ where $K\in \mathbb{R}$. Suppose that $(\ref{PPPPP})$ holds. Let $u(\cdot, t)=P_tf$ be a positive solution to the heat
equation $\partial_t u=Lu$ with $u(\cdot, 0)=f$, $f\in C(M, (0, \infty))$. Define
\begin{eqnarray*}
H_{K}(f, t)=D_K(t)\int_M (f\log f-P_tf\log P_tf )d\mu,
\end{eqnarray*}
where $D_0(t)={1\over t}$ and $D_{K}(t)={2K\over 1-e^{-2Kt}}$ for $K\neq 0$.Then, for all $K\in \mathbb{R}$,
\begin{eqnarray*}
{d\over dt}H_{K}(f, t)\leq 0,\ \ \ \forall t\in (0, T],
\end{eqnarray*}
and for all  $t\in (0, T]$, we have
\begin{eqnarray}
& &{d^2\over dt^2}H_{K}(f, t)+2K\coth(Kt){d\over dt}H_K(f, t)+2D_K(t) \int_M  |\nabla^2
\log P_tf|^2 P_tf d\mu
\nonumber\\
& &\hskip1cm  =
2D_K(t)\int_M \left({1\over 2}{\partial g\over \partial t}+Ric(L)-Kg\right)(\nabla \log P_tf, \nabla \log P_tf)P_tfd\mu.\label{HHH3}
\end{eqnarray}
In partcular, for all $t\in (0, T]$, we have
\begin{eqnarray*}
{d^2\over dt^2}H_K(f, t)+2K\coth(Kt) {d\over dt}H_K(f, t)+2D_K(t)\int_M |\nabla^2\log P_tf|^2P_tfd\mu\leq 0, 
\end{eqnarray*}
and the equality holds for non trivial $f$  in $(0, T]$ if and only if 
\begin{eqnarray*}
{1\over 2}{\partial g\over \partial t}+Ric(L)=Kg, \ \ {\partial\phi\over \partial t}=-\Delta \phi-R+nK, \ \ \ \forall t\in (0, T].\label{KRF4}
\end{eqnarray*}
\end{theorem}

\noindent{\bf Proof}. By Theorem $3.1$ in \cite{LL15}, we have
\begin{eqnarray*}
{d\over dt}\int_M {|\nabla P_tf|^2\over P_tf}d\mu&=&-2\int_M |\nabla^2\log
P_tf|^2P_tfd\mu\\
& &\hskip1cm -2\int_M \left({1\over 2}{\partial g\over \partial
t}+Ric(L)\right)(\nabla \log P_tf, \nabla \log P_tf)P_tfd\mu.
\end{eqnarray*}
Similarly to the proof of $(\ref{HHH2})$ in Theorem \ref{WCDK}  we can prove $(\ref{HHH3})$. 
\hfill $\square$

\medskip
Similarly to Section $3$, we  define the $W$-entropy  for the heat equation of the Witten Laplacian on $(K, \infty)$-super Perelman Ricci flow by the following revised Boltzmann entropy formula
\begin{eqnarray*}
W_K(f, t)=H_K(f, t)+{\sinh(2Kt)\over 2K}{d\over dt}H_K(f, t).
\end{eqnarray*}

Then, we have the following 

\begin{theorem} Let $(M, g(t), \phi(t), t\in [0, T])$ be a  compact $K$-super Perelman  Ricci flow, i.e., 
$${1\over 2}{\partial g\over \partial t}+Ric(L)\geq Kg,$$ where $K\in \mathbb{R}$. Suppose that $(\ref{PPPPP})$ holds.  Then
\begin{eqnarray*}
& &{d\over dt}W_{K}(f, t)+(1+e^{2Kt})\int_M  |\nabla^2
\log P_tf|^2 P_tf d\mu\nonumber\\
& &\hskip1cm = -(1+e^{2Kt})\int_M \left({1\over 2}{\partial g\over \partial t}+Ric(L)-Kg\right)(\nabla\log P_tf, \nabla\log P_tf)P_tfd\mu.\label{WWWW2}
\end{eqnarray*}
In particular, for all $t\in (0, T]$, we have 
\begin{eqnarray*}
{d\over dt}W_{K}(f, t)+(1+e^{2Kt})\int_M |\nabla^2\log P_t f|^2P_tf d\mu\leq 0.
\end{eqnarray*}
Moreover, the following statements are equivalent:\\
(i) for all $t\in (0, T]$, and non constant $f$, 
\begin{eqnarray*}
{d^2\over dt^2}H_K(f, t)+2K\coth(Kt) {d\over dt}H_K(f, t)+2D_K(t)\int_M |\nabla^2\log P_tf|^2P_tfd\mu=0.
\end{eqnarray*}
(ii)  for all $t\in (0, T]$, and non constant $f$, 
\begin{eqnarray*}
{d\over dt}W_{K}(f, t)+(1+e^{2Kt})\int_M |\nabla^2\log P_t f|^2P_tf d\mu=0.
\end{eqnarray*}
(iii)  $(M, g(t), t\in [0, T])$ is a  $K$-Perelman Ricci flow, i.e., 
$${1\over 2}{\partial g\over \partial t}+Ric(L)=Kg, 
\ \ {\partial\phi\over \partial t}=-\Delta \phi-R+nK,\ \ \ \ t\in (0, T].$$
\end{theorem}
{\it Proof}. The  proof is similar the one of Theorem \ref{WCDK1}. \hfill $\square$

\section{A second order differential entropy inequality on $(K, m)$-super Ricci flows}
Note that for any $m>n$, we have
\begin{eqnarray*}
|\nabla^2 P_t f|^2&\geq& {|\Delta \log P_t f||^2\over n}\\
&\geq& {|L\log P_t f|^2\over m}-{|\nabla \log P_t f\cdot \nabla \phi|^2\over m-n},
\end{eqnarray*}
whence 
\begin{eqnarray*}
\int_M |\nabla^2 P_t f|^2 P_tf d\mu\geq {1\over m} \int_M  |L\log P_t f|^2 P_tfd\mu-{1\over m-n}\int_M |\nabla \log P_t f\cdot \nabla \phi|^2P_t fd\mu.
\end{eqnarray*}
By the Cauchy-Schwartz inequality and integration by parts formula, we have
\begin{eqnarray*}
\int_M |L\log P_t f|^2P_t fd\mu\geq \left(\int_M L\log P_t f P_t fd\mu\right)^2=\left(\int_M |\nabla \log P_t f|^2P_tfd\mu\right)^2.
\end{eqnarray*}
Combining this with $(\ref{HHH2})$,  we prove the following  result. 

\begin{theorem} Let $(M, g)$ be a complete Riemannian manifold with the bounded geometry condition, $\phi\in C^4(M)$ with $\nabla\phi\in C_b^3(M)$. Then
\begin{eqnarray*}\label{HHIK1}
& &{d^2\over dt^2}H_K(f, t)+2K\coth(Kt){d\over dt}H_K(f, t)+{2D_K(t)\over m}\left(\int_M {|\nabla P_t|^2\over P_tf} d\mu\right)^2\nonumber\\
& &\hskip1.5cm \leq -2D_K(t)\int_M (Ric_{m, n}(L)-Kg)(\nabla\log P_tf, \nabla\log P_tf)P_tfd\mu.
\end{eqnarray*}
In particular, if $Ric_{m, n}(L)\geq Kg$, we have
\begin{eqnarray*}\label{HHIK1}
{d^2\over dt^2}H_K(f, t)+2K\coth(Kt){d\over dt}H_K(f, t)+{2D_K(t)\over m}\left(\int_M {|\nabla P_t|^2\over P_tf} d\mu\right)^2\leq 0,
\end{eqnarray*}
and the equality holds if and only if  $(M, g, \phi)$ is a $(K, m)$-Ricci soliton, i.e.,  
\begin{eqnarray*}
Ric_{m, n}(L)=Kg.
\end{eqnarray*}
\end{theorem}

Similarly, we have the following second order differential entropy inequality on manifolds with time dependent metrics and potentials.

\begin{theorem} \label{HIK2} Let $(M, g(t), \phi(t), t\in [0, T])$ be a compact manifold equipped with time dependent metrics and potentials which satisfy the constraint equation 
$\partial_t \phi ={1\over 2}{\rm Tr}\partial_t g$. Then 
\begin{eqnarray*}
& &{d^2\over dt^2}H_K(f, t)+2K\coth(Kt){d\over dt}H_K(f, t)+{2D_K(t)\over m}\left(\int_M {|\nabla P_t|^2\over P_tf} d\mu\right)^2\nonumber\\
& &\hskip0.6cm \leq -2D_K(t)\int_M \left({1\over 2}{\partial g\over \partial t}+Ric_{m, n}(L)-Kg\right)(\nabla\log P_tf, \nabla\log P_tf)P_tfd\mu. 
\end{eqnarray*}
In particular, if $(M, g(t), \phi(t), t\in [0, T])$ is a compact manifold equipped with a $(K, m)$-super Perelman Ricci flow, i.e., 
\begin{eqnarray*}
{1\over 2}{\partial g\over \partial t}+Ric_{m, n}(L)\geq Kg,
\end{eqnarray*}
then 
\begin{eqnarray*}
{d^2\over dt^2}H_K(f, t)+2K\coth(Kt){d\over dt}H_K(f, t)+{2D_K(t)\over m}\left(\int_M {|\nabla P_t|^2\over P_tf} d\mu\right)^2\leq 0,
\end{eqnarray*}
and the equality holds if and only if $(M, g(t), \phi(t), t\in [0, T])$ is a $(K, m)$-Perelman Ricci flow, i.e., 
\begin{eqnarray*}
{1\over 2}{\partial g\over \partial t}+Ric_{m, n}(L)=Kg, \ \ \ {\partial \phi\over \partial t}={1\over 2}{\rm Tr}\left({\partial g\over \partial t}\right).
\end{eqnarray*}
\end{theorem}
{\it Proof}. By the same argument as above, we can prove Theorem \ref{HIK2} by using $(\ref{HHH3})$. \hfill $\square$

\section{Probabilistic interpretation of the $W_{m, K}$-entropy}

In this section, we  give a  probabilistic interpretation of the  $W_{m, K}$-entropy for the heat equation of the Witten Laplacian on manifolds with the $CD(K, m)$-condition, where $K\in \mathbb{R}$ and $m\in [n, \infty)\cap \mathbb{N}$. . 

Let $m\in [n, \infty)\cap\mathbb{N}$, $M=\mathbb{R}^m$, $g_0$ the Euclidean metric, $\phi_K(x)=-{K\|x\|^2\over 2}$ and $d\mu_K(x)=e^{K\|x\|^2\over 2}dx$, where $K\in \mathbb{R}$.  Then $\nabla\phi_K(x)=-Kx$, and $\nabla^2\phi_K=-K{\rm Id}_{\mathbb{R}^m}$. We consider the Ornstein-Ulenbeck operator on $\mathbb{R}^m$ given by
$$L=\Delta+Kx\cdot \nabla.$$
By Section $3.3$,  $(\mathbb{R}^m, g_0, \phi_{K})$ is a complete shrinking Ricci soliton, i.e., $Ric(L)=-Kg_0$. The  Ornstein-Ulenbeck diffusion process on $\mathbb{R}^m$ satisfies the Langevin SDE
$$dX_t=\sqrt{2} dW_t+KX_tdt, \ \ \ X_0=x,$$
and is given by the explicit formula below
$$
X_t=e^{Kt}x+\sqrt{2}\int_0^ t e^{K(t-s)}dW_s.$$
Hence 
\begin{eqnarray*}
X_t =e^{Kt}x+\sqrt{e^{2Kt}-1\over K}\xi\ \ \ {\rm in\ law}, 
\end{eqnarray*}
where $\xi$ is a standard $N(0, {\rm Id})$ variable on  $\mathbb{R}^m$. Thus the law of $X_t$ is Gaussian  $N\left(e^{Kt}x, {e^{2Kt}-1\over K}{\rm Id}\right)$, and the heat kernel of $X_t$ with respect to the Lebesgue measure on $\mathbb{R}^m$  is given by
\begin{eqnarray*}
u_{m, K}(x, y, t)=\left({K\over 2\pi (e^{2Kt}-1)}\right)^{m/2}\exp \left(-{K|y-e^{Kt}x|^2\over 2(e^{2Kt}-1)}\right).
\end{eqnarray*}

Fix $x\in \mathbb{R}^m$, and denote $\sigma_K^2={e^{2Kt-1}\over 2K}$. The relative Boltzmann-Shannon entropy of the law of $X_t$ with respect to the Lebesgue measure on $\mathbb{R}^m$  is given by

\begin{eqnarray*}
{\rm Ent}(u_{m, K}(x, y, t)| dy)&=&\int_{\mathbb{R}^m} u_{m, K}(x, y, t)\log u_{m, K}(x, y, t)dy\\
&=&-{m\over 2}\left(1+\log(4\pi \sigma_K^2(t))\right).
\end{eqnarray*}
When   $t\rightarrow 0$, we have 
\begin{eqnarray*}
{\rm Ent}(u_{m, K}(x, y, t) | dy)&=&- {m\over 2}\left[1+ \log\left(4\pi\times {e^{2Kt}-1\over 2K}\right)\right]\\
&=&-{m\over 2}\left(1+\log (4\pi t)+\log(1+Kt+{2K^2t^2\over 3}+{K^3t^3\over 3}+O(t^4)) \right)\\
&=&-{m\over 2}\left(1+\log (4\pi t)+Kt+{K^2t^2\over 6}\right)+O(t^4).
\end{eqnarray*}

Thus, when $t\rightarrow 0^+$, the second term in the definition formula $(\ref{HmK})$ of the $H_{m, K}$-entropy is  asymptotically (with order $O(t^4)$) equivalent to the Boltzmann-Shannon entropy of the heat kernel at time $t$ of the Ornstein-Uhlenbeck operator on $\mathbb{R}^m$ with respect to the Lebesgue measure on $\mathbb{R}^m$.  That is to say,  when $t\rightarrow 0^+$, we have
\begin{eqnarray*}
H_{m, K}(u(t))={\rm Ent}(u_{m, K}(t)|dy)-{\rm Ent}(u(t)|\mu)+O(t^4),
\end{eqnarray*}
while 
\begin{eqnarray*}
W_{m, K}(u(t))={d\over dt}\left(tH_{m, K}(u(t))\right).
\end{eqnarray*}
Moreover,  on complete Riemannian manifolds with the $CD(-K, m)$-condition and on $(-K, m)$-super Perelman Ricci flows, we have 
\begin{eqnarray*}
{d\over dt}W_{m, K}(u(t))\leq 0, \ \ \forall t\in (0, T], \label{WmKt3}
\end{eqnarray*}
and asymptotically when $t\rightarrow 0^+$, we have 
\begin{eqnarray*}
{d^2\over dt^2}(t{\rm Ent}(u(t)|\mu))\geq -{m\over 2t}(1+Kt)^2={d^2\over dt^2}(t{\rm Ent}(u_{m, K}(t)|dy))+O(t^2).
\end{eqnarray*}

\begin{flushleft}
\medskip\noindent

Songzi Li, School of Mathematical Science, Beijing Normal University, No. 19, Xin Jie Kou Wai Da Jie, 100875, China, Email: songzi.li@bnu.edu.cn

\medskip

Xiang-Dong Li, Academy of Mathematics and Systems Science, Chinese
Academy of Sciences, 55, Zhongguancun East Road, Beijing, 100190, China, 
E-mail: xdli@amt.ac.cn
\\
and
\\
School of Mathematical Sciences, University of Chinese Academy of Sciences, Beijing, 100049, China
\end{flushleft}

\end{document}